\documentclass[proceedings,submission%
]{dmtcs}

\usepackage[latin1]{inputenc}
\usepackage{subfigure}
\usepackage{color}
\usepackage{shuffle}
\usepackage{amsfonts}
\usepackage{amssymb}
\usepackage{amscd}
\usepackage{amstext}
\usepackage{amsmath}
\usepackage{rotating}
\usepackage{xy}
\usepackage{mathrsfs}
\usepackage{stmaryrd}
\xyoption{all}

\usepackage[round]{natbib}

\newtheorem{example}{Example(s)}[section]

\newtheorem{theorem}[example]{Theorem}

\newtheorem{proposition}[example]{Proposition}

\newtheorem{remark}[example]{Remark}

\def\hs{\hbox to 3mm{}}
\def\hhs{\hbox to 5cm{}}
\def\ss{\smallskip}

\gdef\stuffle{\;%
  \setlength{\unitlength}{0.0125cm}%
  \begin{picture}(20,10)(220,580)
  \thinlines
  \put(220,592){\line( 0,-1){ 10}}
  \put(220,582){\line( 1, 0){ 20}}
  \put(240,582){\line( 0, 1){ 10}}
  \put(230,592){\line( 0,-1){ 10}}
  \put(225,587){\line( 1, 0){ 10}}
  \end{picture}\;
}

\def\adots{\mathinner{\mkern2mu\raise1pt\hbox{.}
\mkern3mu\raise4pt\hbox{.}\mkern1mu\raise7pt\hbox{.}}}

\def\up#1{\raise 1ex\hbox{\footnotesize#1}}

\def\mref#1{(\ref{#1})}

\def\ra{\rightarrow}

\def\A{\mathcal{A}}

\def\N{{\mathbb N}}

\def\al{\alpha}
\def\be{\beta}
\def\ga{\gamma}

\def\g{{\mathfrak g}}
\def\U{\mathcal{U}}

\def\scal#1#2{\langle #1 | #2 \rangle}
\def\kd#1#2{\delta_{#1#2}}

\def\ncp#1#2{#1\langle #2\rangle}

\def\End{\mathrm{End}}

\author{Matthieu Deneufch\^atel\addressmark{1} \and G\'erard H. E. Duchamp\addressmark{1} \and Vincel Hoang Ngoc Minh\addressmark{1}}
\title[Radford bases and Sch\"utzenberger's Factorizations.]{Radford bases and Sch\"utzenberger's Factorizations.}
\address{\addressmark{1}Laboratoire d'Informatique de Paris Nord, UMR 7030 CNRS, Universit\'e Paris 13, 99, Avenue J.-B. Cl\'ement, 93430 Villetaneuse, France\\
}
\keywords{Sch\"utzenberger's factorizations, Poincar\'e-Birkhoff-Witt basis, Radford basis, Cartier-Quillen-Milnor-Moore Theorem}

\begin{document}
\maketitle
\begin{abstract}

\paragraph{Abstract.}
In this paper, we present Sch\"utzenberger's factorization in different combinatorial contexts and show that its validity is not restricted to these cases but can be extended to every Lie algebra endowed with an ordered basis. We also expose some elements of the relations between the Poincar\'e-Birkhoff-Witt bases of the enveloping algebra and their dual families.

\paragraph{R\'esum\'e.}
Dans cet article, nous pr\'esentons la factorisation de Sch\"utzenberger dans diff\'erents contextes combinatoires avant de montrer que c'est en fait une relation g\'en\'erale valide pour toute alg\`ebre de Lie dot\'ee d'une base ordonn\'ee. Nous pr\'esentons par ailleurs quelques \'el\'ements relatifs aux liens qui unissent les bases de Poincar\'e-Birkhoff-Witt de l'alg\`ebre enveloppante de l'alg\`ebre de Lie consid\'er\'ee et leurs familles duales.

\end{abstract}

\section{Introduction}

On the ground of computation, the Haussdorff group has not received the attention it deserves.
 
\ss
Though its definition is simple (it is the group of group-like elements in a suitable complete 
bialgebra) it seems that, even in the free context (free Lie algebra), this group permits to encompass many desired combinatorial tools for a theory of infinite dimensional (Combinatorial) Lie Groups.

\ss
This is done through Sch\"utzenberger's factorization formula in its resolution-of-unity-like formulation\footnote{Here indexed by the partially commutative monoid $w\in M(X,\theta)$ (\cite{CaFo}).} 
\begin{equation}\label{SF1}
	\sum_{w\in X^*} w\otimes w = \prod_{l \in \rm{Lyn}(X)}^{\searrow} e^{S_l\otimes P_l}
\end{equation}
which provides a beautiful framework for ``local coordinates'' in this infinite dimensional Lie Group. 

\ss
The paper is organized as follows. Section \ref{generalities} is devoted to generalities. In Section \ref{examples}, we present four combinatorial examples : non commutative, partially commutative and commutative free algebras and the case of the stuffle algebra. Finally, in Section \ref{generalsettings}, we give the general theorem and presents the duality between Radford bases and bases of Poincar\'e-Birkhoff-Witt type.\\
 
\section{Generalities}
\label{generalities}
\textbf{Multiindex notation :} If $Y=(y_i)_{i\in I}$ is a totally ordered family in an algebra $\A$ and $\al\in \N^{(I)}$, one defines $Y^{\al}$ by
\begin{equation}
	y_{i_1}^{\al(i_1)}y_{i_2}^{\al(i_2)}\cdots y_{i_k}^{\al(i_k)}
\end{equation}
for every subset $J=\{i_1,i_2\cdots i_k\}\ , \ i_1 > i_2 > \cdots > i_k$, of $I$ which contains the support of $\al$ (it is easily shown that the value of $Y^\al$ does not depend on the choice of $J\supset supp(\al)$).\\
In particular, if $(e_i)_{i \in I}$ denotes the canonical basis of $\N^{(I)}$, one has $Y^{e_i} = y_i$.

\noindent \textbf{Characteristic :}
Throughout the paper, $k$ denotes a field of characteristic $0$.

\section{Combinatorial examples}
\label{examples}
\subsection{Non Commutative case}
Let $X$ be an alphabet (totally) ordered with $<$. We denote by $\rm{Lyn}(X)$ the set of Lyndon words with letters in $X$. The standard factorization (\cite{Reutenauer93}) of $l \in \rm{Lyn}(X)$ is denoted by $\sigma(l) = ( l_1 , l_2 )$ where $l_2$ is the Lyndon proper right-factor of $l$ of maximal length. The standard factorization and the fact that every word $w \in X^*$ can be factorized as a decreasing product of Lyndon words allow us to define a triangular basis $(P_{w})_{w \in X^*}$, of \emph{Poincar\'e-Birkhoff-Witt} type for the free algebra $\ncp{k}{X}$ as follows : 
\begin{equation}
\label{bracketing}
   P_w = \left\{
\begin{aligned}
   w \quad \quad \quad \quad \quad \, \, \, & \, \text{ if }|w|=1 ; \\
\left[P_{l_1} ,P_{l_2} \right]\quad \quad \quad \quad \, & \, \text{ if } w=l \in {\rm Lyn}(X) \text{ and } (l_1 , l_2) = \sigma(l) ; \\ 
P_{l_1}^{\alpha_1} \dots P_{l_n}^{\alpha_n} \quad \quad \quad \, & \, \text{ if } w = l_{i_1}^{\alpha_1} \dots l_{i_k}^{\alpha_k} \text{ with } l_1 > \dots > l_n\ .
\end{aligned}
\right.
\end{equation}
This basis is triangular because one has
$$
P_w = w + \sum_{u > w} \scal{P}{u} u\ .
$$

Because of the multihomogeneity of $(P_{w})_{w \in X^*}$, it is possible to construct a basis $(S_w)_{w \in X^*}$ of $\ncp{k}{X}$ satisfying $\scal{S_u}{P_v} = \kd{u}{v}$, for all $u,v \in X^*$. One can show that (\cite{Reutenauer93}) $S_w$ is given by
\begin{equation}
   S_w = \left\{
\begin{aligned}
   w \quad \quad \quad \quad \, & \, \text{if }|w|=1 ; \\
   x S_u \quad \quad \quad \, \, \, \, \, & \, \text{if }w=xu \text{ and }w \in \text{Lyn}(X) ; \\
   \frac{S_{l_{i_1}}^{\, \shuffle \, \alpha_1} \shuffle \, \dots \shuffle \, S_{l_{i_k}}^{\, \shuffle \, \alpha_k}}{\alpha_1 ! \dots \alpha_k !} \quad \, & \, \text{ if } w = l_{i_1}^{\alpha_1} \dots l_{i_k}^{\alpha_k} \text{ (decreasing factorization)} \ .
\end{aligned}
\right.
\end{equation}

With these notations, the following equality holds :
\begin{equation}
   \sum_{w \in X^*} w \otimes w = \prod_{l \in \rm{Lyn}(X)}^{\searrow} \exp(S_l \otimes P_l)
\end{equation}
where the product in the right-hand side is the shuffle product on the left and the concatenation product on the right.

\subsection{Partially commutative case}

Let $X$ be a set and $\theta \subset X \times X$ a symmetric and antireflexive relation on $X$ (here, antireflexive means that for all $x \in X, \, (x,x) \notin \theta$). We denote by $M(X , \theta)$ the free partially commutative monoid over $X$ (\cite{CaFo}, \cite{vie86}). It is defined by generators and relations by
\begin{equation}
 M(X,\theta)=\langle X , \left\{ ( x y , y x ) \right\}_{(x,y) \in \theta} \rangle_{\rm Mon}\ .
\end{equation}
Let $\ncp{k}{X,\theta}$ denote the partially commutative free algebra over $X$ (\cite{DuchampKrob}), defined by generators and relations by $\langle X, ( xy = yx)_{(x,y) \in \theta} \rangle_{k\text{-alg}}$, and $k \left[ M(X,\theta ) \right]$ the algebra of the partially commutative free monoid. By universal arguments, one can easily see that 
\begin{equation}
 \ncp{k}{X,\theta} \cong k \left[ M(X , \theta ) \right]\ .
\end{equation}
Therefore, it is possible to consider the elements of $k\left[ M(X,\theta ) \right]$ as polynomials over the partially commutative free monoid and set, for all $P \in k \left[ M(X,\theta ) \right]$, 
\begin{equation}
P = \sum_{m \in M(X,\theta)} \scal{P}{m} m\ .
\end{equation}

We are interested in the Hopf algebra structure of $(\ncp{k}{X,\theta} , \mu, 1_{M(X,\theta)} , \Delta , \epsilon , S)$ where ($\mu$ and $1_{M(X,\theta)}$ being straightforward) $\epsilon (P) = \scal{P}{1}$, $\Delta(x) = x \otimes 1 + 1 \otimes x$ and $S(x_1 \dots x_n) = (-1)^{n} x_n \dots x_1$. It is known that the primitive elements of $\ncp{k}{X,\theta}$ are the elements of the free partially commutative Lie algebra ${\mathscr L}_k (X , \theta)$ : ${\rm Prim} (\ncp{k}{X,\theta}) = {\mathscr L}_k (X , \theta)$. 

Moreover, it is possible to generalize Lyndon words to the partially commutative monoids (\cite{lal}) : a partially commutative Lyndon word is a non-empty, primitive (partially commutative) word which is minimal (for the order on $M(X,\theta)$ induced by the lexicographic order on well chosen normal forms (\cite{LalondeKrob})) in its conjugacy class. We denote their set by $\text{Lyn}(X , \theta)$. Krob and Lalonde have generalized the standard factorization of Lyndon words to the partially commutative case :
\begin{proposition}
 Let $w$ belong to $M(X , \theta )$ with length $\geq 2$. Then there exists a unique factorization $w = fn$, called \emph{standard factorization} of $w$, this unique pair being denoted by $\sigma (w) = ( f , n )$, such that
 \begin{enumerate}
  \item $f \neq 1$ ;
  \item $n \in {\rm Lyn}(X, \theta)$ ;
  \item $n$ is minimal among all possible partially commutative Lyndon words that provide a factorization of $w$.
 \end{enumerate}
Moreover, if $l \in {\rm Lyn}(X,\theta)$ with length $\geq 2$ and with $a_i$ as unique initial letter, and if $\sigma(l) = (f,n)$ is the standard factorization of $l$, then $f \in {\rm Lyn}(X , \theta)$ with $a_i$ as unique initial letter and $f < l < n$.
\end{proposition}
These properties allow us to construct a family $(P_l)_{l \in {\rm Lyn}(X,\theta)}$ in the same way as in the commutative case (see Eq. (\ref{bracketing})). One can show (\cite{lal}) that this family forms a basis of ${\mathscr L}_k(X,\theta)$ and that $P_l$ satifies 
\begin{equation}
 P_l = l + \sum_{\genfrac{}{}{0pt}{}{l' > l}{l' \in \text{Lyn}(X,\theta)}} \alpha_{l'} l'.
\end{equation}

Moreover, it is possible to show that each partially commutative word admits a unique nonincreasing factorization in terms of (partially commutative) Lyndon words. Therefore, one can define $P_w, \, w \in M(X,\theta)$ and show that (for example by translating the proof of \cite{Reutenauer93} in the language of partially commutative words)
\begin{equation}
   P_w = w + \sum_{u > w \in M(X,\theta)} \scal{P_w}{u} u\ .
\end{equation}

Finally, the non commutative construction can be extended to the dual family $S_w$ and this allows us to write the following factorization
\begin{equation}
   \sum_{w \in M(X,\theta)} w \otimes w = \prod_{l \in \rm{Lyn}(X , \theta)}^{\searrow} \exp^{S_l \otimes P_l}\ .
\end{equation}

\begin{remark}
   Note that Reutenauer had noticed that the construction of the dual basis is possible in every enveloping algebra (see Theorem 5.3 and Section 5.7 of \cite{Reutenauer93}).
\end{remark}

\subsection{Commutative case}
The commutative case is obtained from the partially commutative setting by choosing $\theta = X \times X - \text{diag}(X)$ (where $\text{diag} (X) = \left\{ (x,x) , \, x \in X \right\}$). Then $\ncp{k}{X,\theta} = k \left[ X \right]$ is the algebra of commutative polynomials and the primitive elements are the homogeneous polynomials of degree 1 : 
\begin{equation}
   \text{Prim}(k \left[ X \right]) = k . X = \left\{ P = \sum_{x \in X} \scal{P}{x} x \right\}\ .
\end{equation}
The set of Lyndon words is $X$. The specialization of equation (\ref{SF1}) yields
\begin{equation}
   \prod_{x \in X} \exp(x \otimes x) = \sum_{\al \in \N^{(X)}} X^\al \otimes X^\al= \sum_{w \in X^{\oplus}} w \otimes w\ .
\end{equation}
Indeed, with $x$ a letter, one has
\begin{equation}
   x^k \, \shuffle \, x = \frac{(k+1)!}{k!} x^{k+1} = (k+1) x^k\ .
\end{equation}
Thus, 
\begin{equation}
   \prod_{x \in X} \sum_{n \geq 0} \frac{x^{\shuffle \, n} \otimes x^n}{n!} = \prod_{x \in X} \sum_{n \geq 0} x^n \otimes x^n
\end{equation}
and one easily recovers the result.

\subsection{Stuffle algebra}
Let $Y = \left\{ y_i \right\}_{i \geq 1}$. We endow $\ncp{k}{Y}$ with the stuffle product given by the following recursion : for all $y_i, \, y_j \in Y$ and for all $u , \, v \in Y^*$,
\begin{equation}
  \left\lbrace 
\begin{aligned}
u \stuffle 1 & = 1 \stuffle u = u ; \\   
y_i u \stuffle y_j v & = y_i ( u \stuffle y_j v ) + y_j ( y_i u \stuffle v ) + y_{i+j} ( u \stuffle v)\ .
\end{aligned}
\right.
\end{equation}
We define on $\ncp{k}{Y}$ a gradation with values in $\N$ given by an integer valued weight function on $Y^*$. It is a morphism of monoids given on the letters by $|y_s| = s$. Thus
\begin{equation}
	|w| = \sum_{k=1}^{\ell(w)}|w\left[ k \right]|\ ,
\end{equation}
the only word of weight $0$ is the empty word and the number of words of weight $n>0$ is $2^{n-1}$. Therefore, $(\ncp{k}{Y} , \stuffle , 1_{Y^*} )$ is graded in finite dimensions. The stuffle product then admits a dual law denoted by $\Delta_{\stuffle}$. It satisfies 
\begin{equation}
   \scal{ u \stuffle v}{w} = \scal{u \otimes v}{\Delta_{\stuffle} (w)}, \, \text{ for all } \, u , \, v , \, w \, \in Y^*\ .
\end{equation}
It is given on the letters by
 \begin{equation}
    \Delta_{\stuffle}(y_s) = y_s \otimes 1  + 1 \otimes y_s + \sum_{s_1+s_2 =s} y_{s_1} \otimes y_{s_2}
 \end{equation}
and one can prove that $\Delta_{\stuffle}$ is a morphism of algebras from $\ncp{k}{Y}$ to $\ncp{k}{Y} \otimes \ncp{k}{Y}$.

Then, $(\ncp{k}{Y} , conc, 1_{Y^*} , \Delta_{\stuffle} , \epsilon)$ is a $\N$-graded cocommutative bialgebra. Thus, it is possible to apply the Cartier-Quillen-Milnor-Moore theorem which ensures that $\ncp{k}{Y}$ is the enveloping algebra of the (Lie-) algebra of its primitive elements (we recall that $\text{char}(k) = 0$) :
\begin{equation}
   \ncp{k}{Y} \equiv {\mathscr U} ( \text{Prim}(\ncp{k}{Y}))\ .
\end{equation}

Let $(B,<)$ be any totally ordered basis of $\text{Prim}(\ncp{k}{Y})$ and 
$(S_\al)_{\al\in \N^{(B)}}$ the dual basis of $(B^\al)_{\al\in \N^{(B)}}$. By the general setting (see below), one has

\begin{equation}
   \prod_{b \in B}^{\searrow} \exp(S_b \otimes b) = \sum_{w \in Y^*} w \otimes w\ .
\end{equation}

\section{General setting.}
\label{generalsettings}
\subsection{From PBW to Radford}

As potential combinatorial applications include : free Lie algebra (noncommutative or with partial commutations as presented above) and finite dimensional Lie algebra where the factorization has only a finite number of terms, we prefer to state the result with its full generality (\textit{i.e.} considering an arbitrary Lie algebra). Let us first give the context. 

\ss
Let $\g$ be a $k$-Lie algebra and $B=(b_i)_{i\in I}$ be an ordered basis of it. The PBW theorem states exactly that $(B^\al)_{\al\in \N^{(I)}}$ is a basis of $\U(\g)$.\\ 
Now, one considers, in  $\U(\g)$, $(S_\al)_{\al\in \N^{(I)}}$, the dual family, i. e. the family of linear forms on $\U$ defined by
\begin{equation}
\scal{S_\al}{B^\beta}=\delta_{\al,\be}\ .
\end{equation}
One has 
\begin{eqnarray}
\label{mult}
S_\al*S_\be \stackrel{(1)}{=} \sum_{\gamma\in \N^{(I)}} \scal{S_\al*S_\be}{B^\gamma}S^\gamma=
\sum_{\gamma\in \N^{(I)}} \scal{S_\al\otimes S_\be}{\Delta(B^\gamma)}^{\otimes 2}S^\gamma=\cr
\sum_{\gamma\in \N^{(I)}} \scal{S_\al\otimes S_\be}
{\sum_{\ga_1+\ga_2=\ga}\frac{\ga !}{\ga_1 !\,\ga_2 !}\, B^{\ga_1}\otimes B^{\ga_2}}^{\otimes 2}S^\gamma=
\frac{(\al+\be) !}{\al !\,\be !}\, S_{\al+\be}
\end{eqnarray}
which shows that the family $T_\al=\al !\, S_\al$ is multiplicative ($T_\al*T_\be=T_{\al+\be}$)\footnote{As this family is (linearly) free, the correspondence $k[I]\ra \U^*$ is an isomorphism onto its image. This image is exactly the space of linear forms that are of finite support {\it on the PBW basis $(B^\al)_{\al\in N^{(I)}}$}.}.

\begin{remark}
 At first, the right-hand-side member of equality $(1)$ of relation $\ref{mult}$ may provide an infinite sum (we do not know whether $S_\al*S_\be \in span_{\gamma \in \N^{(I)}} (S_\gamma)$) as, for a suitable topology, every $\phi \in \U(\g)^*$ reads
 \[
  \phi = \sum_{\gamma \in \N^{(I)}} \scal{\phi}{B^\gamma} S_\gamma\ .
 \]
\end{remark}

\ss
Now, one can see that since the setting of identity \mref{SF1} requires, in general, infinite sums and products, we need to have at our disposal a topology, a convergence criterion or some limiting process. This will be done by endowing 
$\U(\g)$ with the discrete topology and $\End_k(\U(\g))$ with the topology of pointwise convergence. This means that a net $(f_i)_{i\in A}$ ($A$ is a directed set\footnote{A directed (or filtered) set is an ordered set $(A,<)$ such that every pair of elements is bounded above i. e. 
\begin{equation}
(\forall a,b\in A)(\exists c\in A)(a\geq c,\, b\geq c)\ .
\end{equation}
}) converges to $g\in \End_k(\U(\g))$ iff
\begin{equation}
	(\forall b\in \U(\g))(\exists N\in A)(\forall i\geq N)(f_i(b)=g(b))\ .
\end{equation}
This gives the two following derived criteria considering the partial sums and products.\\
A family $(f_i)_{i\in J}$ will be said {\it summable} if the net of partial sums 
$$
S_F=\sum_{j\in F}f_j
$$ 
(for $F$ any finite subset of $J$) converges to some $g\in \End_k(\U(\g))$. This can be formalized as 
\begin{equation}\label{summable}
	(\forall b\in \U(\g))(\exists F \subset_{finite} J)(\forall F')(F\subset F' \subset_{finite} J\Longrightarrow 
	\Big(\sum_{j\in F'}f_j\Big)(b)=g(b))\ .
\end{equation}
Similarly, a family $(f_i)_{i\in J}$ (this time we need that $J$ be totally ordered) will be said {\it mutipliable} (w.r.t. convolution) if the net of partial products 
$$
M_F=\prod^{\ra}_{j\in F}f_j
$$ 
(for $F$ any finite subset of $J$) converges to some $g\in \End_k(\U(\g))$. This will be formalized as 
\begin{equation}\label{multipliable}
	(\forall b\in \U(\g))(\exists F \subset_{finite} J)(\forall F')(F\subset F' \subset_{finite} J\Longrightarrow 
	\Big(\prod^{\ra}_{j\in F'}f_j\Big)(b)=g(b))\ .
\end{equation}
We are now in position to state the general factorization theorem.

\begin{theorem}
Let $k$	be a field of characteristic zero, $\g$	a $k$-Lie algebra, $B=(b_i)_{i\in I}$ be an ordered basis of it and 
$(B^\al)_{\al\in \N^{(I)}}$ be the associated PBW basis. Denoting $(S_\al)_{\al\in \N^{(I)}}$ the dual family of $(B^\al)_{\al\in \N^{(I)}}$ in $\U^*$, one gets the following 
\begin{equation}\label{maim_fact_thm}
	\sum_{\al\in \N^{(I)}} S_\al\otimes B^\al=\prod^{\ra}_{i\in I} \exp\,(S_{e_i}\otimes B^{e_i})
\end{equation}
where $e_i$ denotes the canonical basis of $\N^{(I)}$ (given by $e_i(j)=\delta_{ij}$). 
\end{theorem}

\begin{remark} The two members of \mref{maim_fact_thm} are in fact a resolution of the identity 
through the mapping 
$$
\Phi : V^* \otimes V \rightarrow \End^{\text finite} (V)
$$
which associates to each separated tensor $f\otimes v\in V^* \otimes V$ the endomorphism $\Phi(f\otimes v) : b \mapsto f(b) \cdot v\ .$ This mapping extends by continuity to series and gives $\sum_{\al\in \N^{(I)}} S_\al\otimes B^\al$ as an expression of $Id_\U$.
\end{remark}

\subsection{From Radford to PBW}

In this paragraph, we take the problem the other way round, starting from a family of linear forms (within $\U(\g)$) $(T_\al)_{\al \in \N^{(I)}}$ such that $T_\al \star T_\be = T_{\al + \be}$ (such a family is called a \emph{Radford family}) and that is in duality with some basis $( B^{\left[ \al \right]} )_{\al \in \N^{(I)}}$ of $\U(\g)$ (that is, one has $\displaystyle \scal{T_\al}{B^{\left[ \be \right]}} = \delta_{\al \be}$). Here, the brackets around the multiindex recall that $B^{\left[ \be \right]}$ is not a product as defined in section \ref{generalities}.

\begin{theorem}
 Let $(T_\al)_{\al \in N^{(I)}}$ be a multiplicative basis of $\U(\g)^*$ in duality with some basis $(B^{\left[ \al \right]})_{\al \in N^{(I)}}$ of $\U(\g)$. Then $(B^{\left[ e_i \right]})_{i \in I}$ is a basis of $\g$.
\end{theorem}

\begin{remark}
 This technique is originated from the application of the CQMM Theorem to the stuffle algebra $(\ncp{k}{Y} , conc, 1_{Y^*} , \stuffle , \epsilon )$ with $\text{\rm Prim}(\ncp{k}{Y}) = \g$ and $\ncp{k}{Y} = \U(\g)$. Note that, in that case, $y_p$ is not primitive anymore if $p>1$. Indeed, one can use $\log_*(I)$ which is a projector on the space of primitive
elements $\text{\rm Prim}(k<Y>)$ :
 \begin{equation}
  \begin{aligned}
   \log_*(I)(y_p) & = y_p - \frac{1}{2} \sum_{p_1+p_2=p} y_{p_1}y_{p_2} \\
   & + \frac{1}{3} \sum_{p_1+p_2+p_3=p} y_{p_1}y_{p_2} y_{p_3} \\
   & - \frac{1}{4} \dots
   \end{aligned}
 \end{equation}
For example,
 \begin{equation}
  \begin{aligned}
   \log_*(I)(y_4) & = y_4 - \frac{1}{2} ( y_{1}y_{3} + y_{2}y_{2} + y_{3}y_{1} ) \\
   & + \frac{1}{3} ( y_{1}y_{1} y_{2} + y_{1}y_{2} y_{1} + y_{2}y_{1} y_{1} ) \\
   & - \frac{1}{4} y_1^4\ .
   \end{aligned}
 \end{equation}
\end{remark}

For the stuffle product, the set $Y$ forms a transcendence basis. 

We recall here the recursive definition of the dual product of the stuffle product :
\begin{equation}
 \begin{aligned}
  \Delta_{\stuffle} (y_n w) & = \Delta_{\stuffle}(y_n) \Delta_{\stuffle}(w) ; \\
  \Delta_{\stuffle} (y_n ) & = y_n \otimes 1 + 1 \otimes y_n + \sum_{\genfrac{}{}{0pt}{}{p+q = n}{p,q \geq 1}} y_p \otimes y_q\ .
 \end{aligned}
\end{equation}
One has
\begin{equation*}
  \scal{ u \stuffle v }{ w } = \scal{ u \otimes v }{ \Delta_{\stuffle} w }\ .
\end{equation*}
Thus
\begin{equation}
 \begin{aligned}
y_p u \stuffle y_q v & = \sum_{w \in X^*} \scal{ y_p u \stuffle y_q v }{ w } w \\
  & = \sum_{w \in X^*} \scal{ y_p u \otimes y_q v }{ \Delta_{\stuffle}(w) } w \\
  & = \scal{ y_p u \otimes y_q v }{ \Delta_{\stuffle}(1) } + \sum_{w \in Y^+} \scal{ y_p u \otimes y_q v }{ \Delta_{\stuffle}(w) } w \\
  & = \sum_{\genfrac{}{}{0pt}{}{r \geq 1}{w \in Y^*}} \scal{ y_p u \otimes y_q v }{ \Delta_{\stuffle} (y_r w) } y_r w \\
  & = \sum_{\genfrac{}{}{0pt}{}{r \geq 1}{w \in Y^*}} \scal{ y_p u \otimes y_q v }{ (y_r \otimes 1 + 1 \otimes y_r) \Delta_{\stuffle} (w) } y_r w \\
  & = \sum_{\genfrac{}{}{0pt}{}{r \geq 1}{w \in Y^*}} \scal{ y_p u \otimes y_q v }{ (y_r \otimes 1) \Delta_{\stuffle} (w) } y_r w + \\
  & + \sum_{\genfrac{}{}{0pt}{}{r \geq 1}{w \in Y^*}} \scal{ y_p u \otimes y_q v }{ (1 \otimes y_r) \Delta_{\stuffle} (w) } y_r w \\
  & + \sum_{\genfrac{}{}{0pt}{}{r \geq 1}{w \in Y^*}} \sum_{r_1 + r_2 =r} \scal{ y_p u \otimes y_q v }{ (y_{r_1} \otimes y_{r_2}) \Delta_{\stuffle} (w) } y_r w \\
  & = \scal{ u \otimes y_q v }{ \Delta_{\stuffle} (w) } y_p u + \scal{ y_p u \otimes v }{ \Delta_{\stuffle} ( w ) } y_q v + \scal{ u \otimes v }{ \Delta_{\stuffle} ( w ) } y_{p+q} \\
   & = y_p(u \stuffle y_q v ) + y_q (y_p u \stuffle v ) + y_{p+q} (u \stuffle v)\ .
 \end{aligned}
\end{equation}
This proves that $\left( \text{Lyn}(Y)^{\stuffle \al} \right)_{\al \in \N^{(\text{Lyn}(Y))}}$ is homogeneous with $|y_i|=i$.

Now, one can consider the set of products of the form $\tilde{B}^\al = \left[ (B^{\left[ e_i \right]})_{i \in I} \right]^\al$ and address the question whether $\tilde{B}^\al = B^{\left[ \al \right]}$. This is true in case where $B^\al$ is a PBW basis w.r.t. the chosen order. 

Therefore, when one starts with a Poincar\'e-Birkhoff-Witt basis, the factorization (\ref{maim_fact_thm}) holds. If one starts with a multiplicative family $(T_\al)_{\al \in \N^{(I)}}$, the following identity holds
\begin{equation}
\sum_{\al\in \N^{(I)}} S_\al \otimes \tilde{B}^\al=\prod^{\ra}_{i\in I} \exp\,(S_{e_i} \otimes B^{\left[ e_i \right]})
\end{equation}
(where $T_\alpha = \alpha ! S_\alpha$). But it remains to be proved that the products of the $B^{\left[ e_i \right]}$'s yield the elements $B^{ \left[ \alpha \right]}$ if one wants to have a factorization of the form (\ref{SF1}).

\section{Conclusion}
Though it frequently appears in relation to the free algebra, Sch\"utzenberger's factorization is not a specific property of this structure. On the contrary, it is a very general relation which holds in every enveloping algebra as shown by our theorem. \\
It is also interesting because it underlines the duality between bases of Poincar\'e-Birkhoff-Witt type and Radford bases. We have not yet fully investigated the relations between these structures and hope to shed more light on this subject by computing more combinatorial examples.

\acknowledgements
\label{sec:ack}
The authors wish to acknowledge support from Agence Nationale de la Recherche (Paris, France) under Program No. ANR-08-BLAN-0243-2 as well as support from ``Projet interne au LIPN'' ``Polyz\^eta functions''.

\bibliographystyle{abbrvnat}
\bibliography{FPSAC}

\begin{thebibliography}{6}
\providecommand{\natexlab}[1]{#1}
\providecommand{\url}[1]{\texttt{#1}}
\expandafter\ifx\csname urlstyle\endcsname\relax
  \providecommand{\doi}[1]{doi: #1}\else
  \providecommand{\doi}{doi: \begingroup \urlstyle{rm}\Url}\fi

\bibitem[Cartier and Foata(1969)]{CaFo}
P.~Cartier and D.~Foata.
\newblock \emph{Probl{\`e}mes combinatoires de commutation et
  r{\'e}arrangements}, volume~85 of \emph{LNM}.
\newblock Springer-Verlag, Berlin, 1969.

\bibitem[Duchamp and Krob(1992)]{DuchampKrob}
G.~H.~E. Duchamp and D.~Krob.
\newblock The lower central series of the free partially commutative group.
\newblock \emph{Semigroup Forum}, 45\penalty0 (3):\penalty0 385--394, 1992.

\bibitem[Krob and Lalonde(1993)]{LalondeKrob}
D.~Krob and P.~Lalonde.
\newblock Partially commutative {Lyndon} words.
\newblock In \emph{10th Annual Symposium on Theoretical Aspects of Computer
  Science}, volume 665 of \emph{LNCS}, pages 237--246, W{\"u}rzburg, Germany,
  Feb 1993. Springer.

\bibitem[Lalonde(1993)]{lal}
P.~Lalonde.
\newblock {Bases de Lyndon des algèbres de Lie libres partiellement
  commutatives}.
\newblock \emph{Theoretical Computer Science}, 117\penalty0 (1-2):\penalty0 217
  -- 226, 1993.
\newblock ISSN 0304-3975.
\newblock \doi{10.1016/0304-3975(93)90315-K}.
\newblock URL
  \url{http://www.sciencedirect.com/science/article/pii/030439759390315K}.

\bibitem[Reutenauer(1993)]{Reutenauer93}
C.~Reutenauer.
\newblock \emph{Free Lie Algebras}.
\newblock Number~7 in London Math. Soc. Monogr. (N.S.). Oxford University
  Press, 1993.

\bibitem[Viennot(1986)]{vie86}
G.~X. Viennot.
\newblock Heaps of pieces {I}: Basic definitions and combinatorial lemmas.
\newblock In G.~Labelle et~al., editors, \emph{Proceedings Combinatoire
  {\'e}numerative, Montr{\'e}al, Qu{\'e}bec (Canada) 1985}, number 1234 in
  Lecture Notes in Mathematics, pages 321--350, Heidelberg, 1986.
  Springer-Verlag.

\end{thebibliography}
\label{sec:biblio}

\end{document}